\documentclass[3p,preprint]{elsarticle}
\usepackage{graphicx}
\usepackage{amsfonts,amsmath,amssymb,amsbsy,enumerate}
\usepackage{fancyhdr}

\makeatletter
\def\ps@pprintTitle{%
  \let\@oddhead\@empty
  \let\@evenhead\@empty
  \let\@oddfoot\@empty
  \let\@evenfoot\@oddfoot
}
\makeatother
\usepackage[normalem]{ulem}
\usepackage[mathscr]{euscript}
\usepackage{pstricks}
\usepackage{multirow}
\usepackage{verbatim}
\usepackage{caption}
\usepackage{slashbox}
\usepackage{subfigure,psfrag}
\usepackage{algorithm}
\usepackage[noend]{algpseudocode}
\usepackage{lineno}
\graphicspath{{./Figuras/}} 

\let\And\undefined
\let\Or\undefined

\usepackage{tikz}
\usetikzlibrary{positioning,chains,fit,shapes,calc}
\definecolor{myblue}{RGB}{80,80,160}
\definecolor{mygreen}{RGB}{80,160,80}
\definecolor{myred}{RGB}{255,0,0}
\definecolor{mybrown}{RGB}{165,42,42}
\usepackage{ulem}

\algrenewcommand\algorithmicthen{:}
\algrenewcommand\algorithmicdo{:}
\algnewcommand\And{\textbf{ and }}
\algnewcommand\Or{\textbf{ or }}
\algnewcommand\New{\textbf{ new }}
\algnewcommand\To{\textbf{ to }}
\algnewcommand\DownTo{\textbf{ down to }}
\algnewcommand\Break{\textbf{break}}
\algnewcommand\Null{\textbf{null}}
\algnewcommand\True{\textbf{true}}
\algdef{SE}[SUBALG]{Indent}{EndIndent}{}{\algorithmicend\ }%
\algtext*{Indent}
\algtext*{EndIndent}

\usepackage{amsthm}
\newtheorem{theorem}{Theorem}

\newtheorem{proposition}{Proposition}
\newtheorem{question}{Question}

\newtheorem{conjecture}{Conjecture}
\newcommand{\etal}{\textit{et~al.}}

\newtheoremstyle{case}{}{}{}{}{\bfseries}{:}{ }{}
\theoremstyle{case}
\newtheorem{case}{Case}
\newtheorem{subcase}{Case}
\numberwithin{subcase}{case}

\begin{document}
	
	\title{On Tuza's conjecture in co-chain graphs} 
	\author[1]{Luis Chahua\fnref{fn2}\corref{cor1}}
	\ead{luis.chahua@utec.edu.pe}
	\author[1]{Juan Guti\'errez\fnref{fn2,fn3}}
	\ead{jgutierreza@utec.edu.pe}

    \cortext[cor1]{Corresponding author}	
     
    \fntext[fn2]{This research has been partially supported by Fondo Semilla UTEC, Problemas estructurales en Teoría~de~Grafos.  871075-2022.}
     	
	\fntext[fn3]{This research has been partially supported by Movilizaciones para Investigación AmSud, PLANarity and distance
IN Graph theory. E070-2021-01-Nro.6997.}

	\address[1]{Departamento de Ciencia de la Computaci\'on, Universidad de Ingenier\'ia y Tecnolog\'ia (UTEC),Lima, Per\'u}

\begin{abstract}
In 1981, Tuza conjectured that the cardinality of a minimum set of
	edges that intersects every triangle of a graph is at most twice the cardinality
	of a maximum set of edge-disjoint triangles.
	This conjecture have been proved for several
	important graph classes, as planar graphs, tripartite graphs, among others. However, it remains open on other important classes of graphs, as chordal graphs. Furthermore, it remains open for main subclasses of chordal graphs, as split graphs and interval graphs.
  In this paper, we show that Tuza's conjecture is valid for co-chain graphs with even number of vertices in both sides of the partition, a known subclass of interval graphs.
\end{abstract}


\maketitle
\thispagestyle{fancy}
\section{Introduction}\label{sec:intro}

In this paper, all graphs considered are simple and the notation and terminology are standard.
A \textit{triangle hitting}
of a graph~$G$ is a set of edges of~$G$ whose removal results in a triangle-free graph; and a \textit{triangle packing} of~$G$ is a set of
edge-disjoint triangles of~$G$. We denote by~$\tau(G)$ (resp.~$\nu(G)$) the cardinality of a minimum triangle hitting (resp. maximum
triangle packing) of~$G$. Tuza posed the following conjecture.

\begin{conjecture}[\cite{Tuza81}]\label{conjecture:tuza}
  For every graph~\(G\), we have~\(\tau(G)\leq 2\,\nu(G)\).
\end{conjecture}%

This conjecture was verified for many classes of graphs~\cite{Botler2020,Cui09,Haxell99,Haxell12,tuza1990}.
However, Tuza's conjecture remains open for several important graph classes, as chordal graphs.
Botler \etal~showed that 
Tuza's conjecture is valid for $K_8$-free chordal graphs \cite[Corollary 3.6]{Botler2020}. But the conjecture is still open for several other important subclasses of chordal graphs, as split graphs and interval graphs.

In this direction, Bonamy \etal~\cite{Bonamy2021} verified this conjecture for threshold graphs, that is, graphs that are both split graphs and cographs.
Another important subclass of chordal graphs are interval graphs.
Several algorithmic and structural problems have been studied in co-chain graphs \cite{Boyaci2018,Kijima2012}, a known subclass of unit-interval graphs \cite[Theorem 9.1.2]{Brandstadt1999}. 
Bonamy \etal~showed Conjecture 1 is valid for co-chain graphs with both sides of the same size divisible
by 4. In this paper, we omit the restriction of
equality between the sizes of the partition.
Moreover, we show that Tuza's conjecture is valid for every co-chain graph when both sides are divisible by 2.

\section{Preliminaries} \label{sec:bcg}

Given a graph $G$ and a vertex $v$ in $G$, we denote by 
$N(v)$ the set of neighbors of $v$.
If $X\subseteq V(G)$, then we denote by $G[X]$ the subgraph induced by $X$.
A \textit{co-bipartite} graph is a graph whose complement is bipartite.
That is, if $G$ is a co-bipartite graph, then there is a partition $\{K_1,K_2\}$ of $V(G)$
such that $G[K_1]$ and $G[K_2]$ are cliques.
We denote a co-bipartite graph $G$ with such partition as $G[K_1,K_2]$.
A \textit{co-chain} graph is a co-bipartite graph
in which the neighborhoods of the vertices in each side can be linearly ordered with respect to inclusion.
That is, if $G[K_1,K_2]$ is a co-chain graph,
we can rename the vertices of $K_1$ as
$c_1,c_2,\ldots,c_{|K_1|}$ and
the vertices of $K_2$ as
$d_1,d_2,\ldots,d_{|K_2|}$
such that $N(c_{i+1}) \subseteq N(c_i)$ for
$i=1,2,\ldots, |K_1|-1$ and
$N(d_{i}) \subseteq N(d_{i+1})$ for
$i=1,2,\ldots, |K_2|-1$.

Given a co-chain graph $G[K_1,K_2]$,
Bonamy \etal~\cite{Bonamy2021} showed that 
if $|K_1|=$~$|K_2|$ and $|K_1|$ is divisible by four, then
$\tau(G) \leq 2\nu (G)$ (they named this class of co-chain graphs, \textit{even balanced co-chain graphs}).
In this paper we show that in fact Tuza's conjecture is valid
for any co-chain graph with both sizes multiple of two. Note that we also omit the restriction of equality between $|K_1|$ and $|K_2|$.

\begin{theorem}
If $G[K_1,K_2]$ is a co-chain graph where $|K_1|$
and $|K_2|$ are even, then $\tau(G) \leq 2\nu (G)$.
\end{theorem}

Given a graph $G$ and two sets of vertices $S$
and $K$, we let $p_G(S,K)$ be a maximum triangle packing in $G[K\cup S]$ among those that do not consider any edge in $G[S]$ and let $p_G(K)$ be a maximum triangle packing in $G[K]$.
If the context is clear, we denote 
them by $p(S,K)$ and $p(K)$, respectively.
We say that two set of vertices in a graph
are \textit{complete to each other} if
they are disjoint and
every vertex in one set is adjacent to
any vertex of the other set.

For our proof, we will use two known results.

\begin{proposition}[see {\cite[Corollary 5]{Bonamy2021}}]\label{prop:pSK}
Let $K$ be a clique and $S$ be an independent set such that they are complete to each other in a graph $G$.
Then 
$|p(S,K)| \geq \frac{| K | - 1}{2}.\min\{| S |, | K |\}$. Moreover, if $|K|$ is even, then
$|p(S,K)| \geq \frac{|K |}{2}.\min\{| S|, |K | - 1\}$.
\end{proposition}

\begin{proposition}[{\cite[
Theorem 2]{Feder2012}}]
\label{prop:pK}
$|p(K_n)|=\frac{1}{3}({{n}\choose{2}} - k)$
where
$$
k=
\begin{cases}
0 & \text{: }n \text{ MOD } 6 \in \{1,3\}\\
4 & \text{: }n \text{ MOD } 6 =5\\
 \frac{n}{2} & \text{: }n \text{ MOD } 6 \in \{0,2\}\\
\frac{n}{2} + 1 & \text{: } n \text{ MOD } 6 =4.
\end{cases}
$$
\end{proposition}

Observe that, by Proposition \ref{prop:pK}, $|p(K_n)| \geq \frac{1}{3}({{n}\choose{2}} - \frac{n}{2} - \frac{3}{2})$
for every $n\geq 1$, $|p(K_n)| \geq \frac{1}{3}({{n}\choose{2}} - 4)$ for every odd $n$, and
$|p(K_n)| \geq \frac{1}{3}({{n}\choose{2}} - \frac{n}{2} - 1)$
for every $n \neq 5$.


\section{Proof of Theorem 1}

For the rest of this section, we fix a co-chain graph~$G[K_\ell, K_m]$.
Also we define the vertices of $K_\ell$ as $c_1,c_2,\ldots,c_{|K_\ell|}$ and
the vertices of $K_m$ as
$d_1,d_2,\ldots,d_{|K_m|}$
such that $N(c_{i+1}) \subseteq N(c_i)$ for
$i=1,2,\ldots, |K_\ell|-1$ and
$N(d_{i}) \subseteq N(d_{i+1})$ for
$i=1,2,\ldots, |K_m|-1$.
Recall that~$K_{\ell}$ and~$K_{m}$ have size divisible by 2. We abuse notation and let $\ell:=|K_\ell|/2$,
and $m:=|K_m|/2$.

\begin{figure}[h!]
\center
\centering
    \captionsetup{justification=centering}
    \resizebox{0.43\textwidth}{!}{
\begin{tikzpicture}[thick,
  every node/.style={draw,circle},
  fsnode/.style={fill=myblue},
  ssnode/.style={fill=mygreen},
  fsnode1/.style={fill=myred},
  fsnode2/.style={fill=mybrown},
  t1/.style={rectangle,draw,inner sep=19pt,text width=1cm, rounded corners=2mm},
  t2/.style={rectangle,draw,inner sep=8pt,text width=1cm, rounded corners=2mm},
  t3/.style={rectangle,draw,inner sep=22pt,text width=1cm, rounded corners=2mm},
  t4/.style={rectangle,draw,inner sep=8pt,text width=1cm, rounded corners=2mm},
  t5/.style={rectangle,draw,inner sep=4.5pt,text width=1cm, rounded corners=2mm},
  -,shorten >= 1pt,shorten <= 1pt
]

\begin{scope}[start chain=going below,node distance=4mm]
\foreach \i in {1,2,...,4}
  \node[fsnode,on chain] (l_\i) {};
\end{scope}

\begin{scope}[xshift=4cm,yshift=1.5cm,start chain=going below,node distance=4mm]
\foreach \i in {1,2,...,8}
  \node[fsnode,on chain] (m_\i) {};
\end{scope}
\node [t1,dashed,myblue,fit=(l_1) (l_4),label={[label distance=0.2cm]above:$\textcolor{myblue}{K_\ell}$}] {};
\node [t2,dashed,myred,fit=(l_1) (l_3),label={[label distance=1.3cm,anchor=north]125:$\textcolor{myred}{X_\ell}$}] {};
\node [t5,dashed,mygreen,fit=(l_1) (l_2),label={[label distance=0.4cm]left:$\textcolor{mygreen}{K_\ell^{top}}$}] {};
\node [t5,dashed,mygreen,fit=(l_3) (l_4),label={[label distance=0.4cm]left:$\textcolor{mygreen}{K_\ell^{bot}}$}] {};
\node [t3,dashed,myblue,fit=(m_1) (m_8),label=above:$\textcolor{myblue}{K_m}$] {};
\node [t4,dashed,myred,fit=(m_4) (m_8),label={[label distance=1.65cm, anchor=north]58:$\textcolor{myred}{X_m}$}] {};
\node [t5,dashed,mygreen,fit=(m_5) (m_8),label={[label distance=0.5cm]right:$\textcolor{mygreen}{K_m^{bot}}$}] {};
\node [t5,dashed,mygreen,fit=(m_1) (m_4),label={[label distance=0.5cm]right:$\textcolor{mygreen}{K_m^{top}}$}] {};
\draw (l_1) -- (m_5);a
\draw (l_1) -- (m_6);a
\draw (l_1) -- (m_7);a
\draw (l_1) -- (m_8);a
\draw (l_2) -- (m_5);a
\draw (l_2) -- (m_6);a
\draw (l_2) -- (m_7);a
\draw (l_2) -- (m_8);a
\draw (l_1) -- (m_4);a
\draw (l_2) -- (m_4);a
\draw (l_1) -- (m_3);a
\draw (l_1) -- (m_2);a
\draw (l_1) -- (m_1);a
\draw (l_3) -- (m_5);a
\draw (l_3) -- (m_6);a
\draw (l_3) -- (m_7);a
\draw (l_3) -- (m_8);a
\draw (l_4) -- (m_7);a
\draw (l_4) -- (m_8);a
\end{tikzpicture}}
\caption{A co-chain graph, with $\ell=2,m=4,x_\ell=3$, and $x_m=5$. For simplicity, the edges inside $K_{\ell}$ and $K_{m}$ are not shown.}
\label{fig:cochainXlXmKltKlbKmtKmb}
\end{figure}
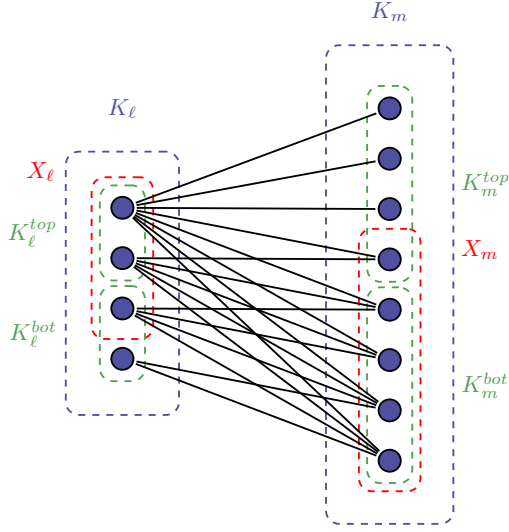

As in \cite{Bonamy2021}, we let~$K^{top}_{\ell}$, $K^{bot}_{\ell}$ be the top and the bottom half of~$K_{\ell}$, respectively, and similarly~$K^{top}_{m}$ and~$K^{bot}_{m}$ be the top and the bottom half of~$K_m$. Let~$X_{\ell} \subseteq K_{\ell}$,~$X_{m} \subseteq K_{m}$ be the sets defined as follows: $c \in X_{\ell}$ if $K^{bot}_{m} \subseteq N(c)$, and~$d \in X_{m}$ if $K^{top}_{\ell} \subseteq N(d)$. We set $x_\ell = | X_{\ell} |$ and $x_m = |X_{m} |$ (Figure \ref{fig:cochainXlXmKltKlbKmtKmb}). As $G$ is a co-chain graph, $x_\ell \geq \ell$ implies that the set $K^{top}_\ell$ is complete to $K^{bot}_{m}$. Consequently, $x_{m} \geq m$. Similarly, $x_{m} \geq m$ implies that $x_{\ell} \geq \ell$. Therefore, $x_{\ell} \geq \ell$ if and only if  $x_{m} \geq m$.
The rest of the proof is divided in two cases, whether $x_{\ell} \geq \ell$ or not.

\subsection{The case $x_\ell \geq \ell$.}

We define a triangle hitting $T_1$ by the union of all edges within $K^{top}_{\ell}$, $K^{bot}_{\ell}$, $K^{top}_{m}$, $K^{bot}_{m}$ as well as all edges between $K^{top}_{\ell}$ and $K^{bot}_{m}$, and between $K^{bot}_{\ell}$ and $K^{top}_{m}$. 
We have
\begin{equation*}
\lvert T_1 \rvert \leq \ell m + (x_\ell  - \ell)(x_m - m)+
2\binom{m}{2}+ 2\binom{\ell}{2}.
\end{equation*}

We will solve first the case $\ell=1$. Note that, if $m\leq 3$, then we are done by Theorem 1.2 of \cite{Puleo15}. So, let us assume that $m\geq 4$. Now, since $2\ell \geq x_{\ell} \geq \ell$, we have $x_{\ell} \in \{1,2\}$.

\begin{figure}[!ht]
  \centering
  \begin{tikzpicture}[scale=0.5]
    \draw[rounded corners=5pt] (0,0) rectangle (4,7);
    \draw (0,2.5) -- (4,2.5);
    
    \draw[rounded corners=5pt] (6,0) rectangle (10,7);
    \draw (6,3.5) -- (10,3.5);
    
    \draw[rounded corners=5pt, red, fill=red!20] (0.2,2.7) rectangle (3.8,6.8);
    \node[font=\footnotesize] at (2,4.75) {$X_{\ell}$};
    \draw[rounded corners=5pt, green] (0.2,0.2) rectangle (3.8,2.3);
    \node[font=\footnotesize] at (2,1.25) {$K_{\ell} \setminus X_{\ell}$};
    
    \draw[rounded corners=5pt, blue, fill=blue!20] (6.2,3.7) rectangle (9.8,6.8);
    \node[font=\footnotesize]  at (8,5.25) {$K_m^{top}$};
    
    \draw[rounded corners=5pt, red, fill=red!20] (6.2,0.2) rectangle (9.8,3.3);
    \node[font=\footnotesize] at (8,2.75) {};
    
    \begin{scope}
      \clip (4,7) -- (4,2.5) -- (6,0) -- (6,3.5);
      \fill[red!10] (0,0) rectangle (10,7);
    \end{scope}
    
    \draw[blue, line width=0.5pt] (7.6,2.2) -- (6.6,4.2);
    \draw[blue, line width=0.5pt] (7.7,2.25) -- (7.1,4.2);
    \draw[blue, line width=0.5pt] (8,2.3) -- (7.8,4.2);
    \draw[blue, line width=0.5pt] (8.1,2.25) -- (8.3,4.2);
    \draw[blue, line width=0.5pt] (8.3,2.2) -- (9,4.2);
    \draw[blue, line width=0.5pt] (8.4,2.15) -- (9.5,4.2);
    
    \node[circle, fill=black, inner sep=0.8pt] (P1) at (6.6,4.2) {};
    \node[circle, fill=black, inner sep=0.8pt] (Q1) at (7.1,4.2) {};
    \node[circle, fill=black, inner sep=0.8pt] (P2) at (7.8,4.2) {};
    \node[circle, fill=black, inner sep=0.8pt] (Q2) at (8.3,4.2) {};
    \node[circle, fill=black, inner sep=0.8pt] (P3) at (9,4.2) {};
    \node[circle, fill=black, inner sep=0.8pt] (Q3) at (9.5,4.2) {};

    \draw[blue, line width=0.5pt] (P1) -- (Q1);
    \draw[blue, line width=0.5pt] (P2) -- (Q2);
    \draw[blue, line width=0.5pt] (P3) -- (Q3);
    
    \draw[blue, line width=0.5pt] (8,1.5) circle (0.8);
    \node[font=\footnotesize] at (8,1.5) {$K_m^{bot}$};
  \end{tikzpicture}
  \caption{Packing $P_2 := p(K_{m}^{bot}\cup X_{\ell}) \cup p(K_{m}^{bot},K_{m}^{top}).$ }
  \label{P2}
\end{figure}
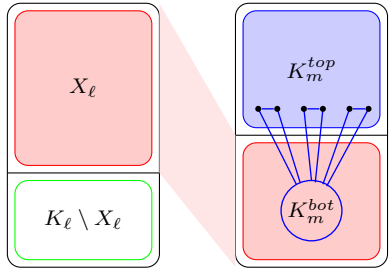

If $x_{\ell}=1$, then $|T_1| \leq m + 2\binom{m}{2} = m^{2}$.
Let~$P_1 := p(K_m).$
By Proposition \ref{prop:pK}, we have
\begin{eqnarray}
2|P_1|- |T_1|
 &\geq& \frac{2}{3}(\binom{2m}{2} - \frac{2m}{2} - 1) -m^{2} \nonumber \\
 &=& \frac{1}{3}(m^{2} -4m-2) \nonumber \\
 &\geq& - \frac{2}{3}. \nonumber
\end{eqnarray} 
Since $2|P_1|- |T_1|$ is integer, we conclude that $2|P_1|- |T_1| \geq 0$. Now, if $x_{\ell}=2$, then $|T_1| \leq 2\binom{m}{2} + x_m \leq~m^{2}+m$. We let
$P_2 := p(K_{m}^{bot}\cup X_{\ell}) \cup p(K_{m}^{bot},K_{m}^{top})$ (Figure \ref{P2}).
By Propositions \ref{prop:pSK} and \ref{prop:pK}, we have
\begin{eqnarray}
2|P_2|- |T_1|
 &\geq& \frac{2}{3}(\binom{m+2}{2} - \frac{m+2}{2} - 1) +m(m-1)-m^{2}-m \nonumber \\
 &=& \frac{1}{3}(m^{2} -4m-2) \nonumber \\
 &\geq& - \frac{2}{3}. \nonumber
\end{eqnarray}
Since $2|P_2|- |T_1|$ is integer, we conclude that $2|P_2|- |T_1| \geq 0$.

\begin{figure}[!h]
  \centering
  \begin{tikzpicture}[scale=0.5]
    \draw[rounded corners=5pt] (0,0) rectangle (4,7);
    \draw (0,3.5) -- (4,3.5);
    \draw[rounded corners=5pt] (6,0) rectangle (10,7);
    \draw (6,3.5) -- (10,3.5);
    
    \draw[rounded corners=5pt, blue, fill=blue!20] (0.2,3.7) rectangle (3.8,6.8);
    \node[font=\footnotesize] at (2,5.38) {$K_{\ell}^{top}$};
    \draw[blue, line width=0.5pt] (2,5.25) circle (0.8);
    \draw[rounded corners=5pt, blue, fill=blue!20] (0.2,0.2) rectangle (3.8,3.3);
    \node[font=\footnotesize] at (2,1.75) {$K_{\ell}^{bot}$};
    
    \draw[rounded corners=5pt, red] (0.1,2.5) rectangle (3.9,6.9);
    \node[font=\footnotesize, text=red] at (-0.5,5.38) {$X_{\ell}$};

    \draw[rounded corners=5pt, green] (6.2,3.7) rectangle (9.8,5.3);
     \node[font=\footnotesize] at (8,4.5) {\resizebox{!}{0.18cm}{$X_m \setminus K_m^{bot}$}};
    \draw[rounded corners=5pt, fill=blue!20] (6.2,0.2) rectangle (9.8,3.3);
    \node at (8,2.75) {};

    \draw[blue, line width=0.5pt] (7.6,2.2) -- (6.4,5.1);
    \draw[blue, line width=0.5pt] (7.7,2.25) -- (6.4,5.7);
    \draw[blue, line width=0.5pt] (8.3,2.2) -- (9.5,5.7);
    \draw[blue, line width=0.5pt] (8.4,2.15) -- (9.5,5.1);
    
    \draw[blue, line width=0.5pt] (1.5,4.6) -- (0.6,2.8);
    \draw[blue, line width=0.5pt] (1.6,4.55) -- (0.6,2.2);
    \draw[blue, line width=0.5pt] (1.8,4.45) -- (2.2,2.2);
    \draw[blue, line width=0.5pt] (1.9,4.45) -- (2.2,2.8);
    \draw[blue, line width=0.5pt] (2.1,4.5) -- (3.4,2.2);
    \draw[blue, line width=0.5pt] (2.2,4.53) -- (3.4,2.8);
    
    \draw[red, line width=0.5pt] (3.9,3.6) -- (6.5,0.5);
    \draw[red, line width=0.5pt] (3.9,3.6) -- (6.5,1.0);
    \draw[red, line width=0.5pt] (3.9,4.0) -- (6.5,1.5);
    \draw[red, line width=0.5pt] (3.9,4.0) -- (6.5,2.0);
    \draw[red, line width=0.5pt] (3.9,4.4) -- (6.5,2.5);
    \draw[red, line width=0.5pt] (3.9,4.4) -- (6.5,3.0);
    
    \draw[blue, line width=0.5pt] (6.3,5.2) -- (3.5,6.6);
    \draw[blue, line width=0.5pt] (6.3,5.2) -- (3.5,6.2);
    \draw[blue, line width=0.5pt] (6.2,4.6) -- (3.5,5.9);
    \draw[blue, line width=0.5pt] (6.2,4.6) -- (3.5,5.5);
    \draw[blue, line width=0.5pt] (6.3,3.9) -- (3.5,5.2);
    \draw[blue, line width=0.5pt] (6.3,3.9) -- (3.5,4.8);

    \node[circle, fill=black, inner sep=0.8pt] (P1) at (6.4,5.1) {};
    \node[circle, fill=black, inner sep=0.8pt] (Q1) at (6.4,5.7) {};
    \node[circle, fill=black, inner sep=0.8pt] (P3) at (9.5,5.7) {};
    \node[circle, fill=black, inner sep=0.8pt] (Q3) at (9.5,5.1) {};
    
    \node[circle, fill=black, inner sep=0.8pt] (R1) at (0.6,2.8) {};
    \node[circle, fill=black, inner sep=0.8pt] (S1) at (0.6,2.2) {};
    \node[circle, fill=black, inner sep=0.8pt] (R2) at (2.2,2.8) {};
    \node[circle, fill=black, inner sep=0.8pt] (S2) at (2.2,2.2) {};
    \node[circle, fill=black, inner sep=0.8pt] (R3) at (3.4,2.8) {};
    \node[circle, fill=black, inner sep=0.8pt] (S3) at (3.4,2.2) {};
    
    \node[circle, fill=black, inner sep=0.8pt] (T1) at (3.5,6.6) {};
    \node[circle, fill=black, inner sep=0.8pt] (U1) at (3.5,6.2) {};
    \node[circle, fill=black, inner sep=0.8pt] (T2) at (3.5,5.9) {};
    \node[circle, fill=black, inner sep=0.8pt] (U2) at (3.5,5.5) {};
    \node[circle, fill=black, inner sep=0.8pt] (T3) at (3.5,5.2) {};
    \node[circle, fill=black, inner sep=0.8pt] (U3) at (3.5,4.8) {};
    
    \node[circle, fill=black, inner sep=0.8pt] (V1) at (6.5,3.0) {};
    \node[circle, fill=black, inner sep=0.8pt] (W1) at (6.5,2.5) {};
    \node[circle, fill=black, inner sep=0.8pt] (V2) at (6.5,2.0) {};
    \node[circle, fill=black, inner sep=0.8pt] (W2) at (6.5,1.5) {};
    \node[circle, fill=black, inner sep=0.8pt] (V3) at (6.5,1.0) {};
    \node[circle, fill=black, inner sep=0.8pt] (W3) at (6.5,0.5) {};     
    
    \draw[blue, line width=0.5pt] (P1) -- (Q1);
    \draw[blue, line width=0.5pt] (P3) -- (Q3);
    
    \draw[blue, line width=0.5pt] (R1) -- (S1);
    \draw[blue, line width=0.5pt] (R2) -- (S2);
    \draw[blue, line width=0.5pt] (R3) -- (S3);
    
    \draw[blue, line width=0.5pt] (T1) -- (U1);
    \draw[blue, line width=0.5pt] (T2) -- (U2);
    \draw[blue, line width=0.5pt] (T3) -- (U3);
    
    \draw[red, line width=0.5pt] (V1) -- (W1);
    \draw[red, line width=0.5pt] (V2) -- (W2);
    \draw[red, line width=0.5pt] (V3) -- (W3);
    
    \draw[blue, line width=0.5pt] (8,1.5) circle (0.8);
    \node[font=\footnotesize][black] at (8,1.5) {$K_m^{bot}$};
    \node at (8,4.85) {};
   \draw (6,5.5) -- (10,5.5);

  \end{tikzpicture}
  \caption{Packing $P_3 := p(X_\ell,K_{m}^{bot}) \cup p(X_m \setminus K_{m}^{bot}, K_{\ell}^{top}) \cup p(K_{m}^{bot}, K_{m}^{top}) \cup p(K_{\ell}^{top}, K_{\ell}^{bot}).$ }
  \label{P3}
\end{figure}
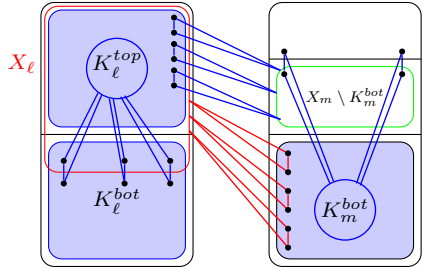

So, we can assume that $\ell \geq 2$. Analogously, we may assume that $m \geq 2$.
For the rest of the proof in this case, we may assume, without loss of generality, that $\ell \leq m$. We divide the rest of the proof in two cases: either $x_\ell \leq m$ or not.
As we will use Propositions \ref{prop:pSK} and  \ref{prop:pK} as before several times, we
will do it without stating them throughout the rest of the paper. \\\\
\newpage
We define a triangle packing $P_3$ as follows (Figure \ref{P3}).
$$P_3 := p(X_\ell,K_{m}^{bot}) \cup p(X_m \setminus K_{m}^{bot}, K_{\ell}^{top}) \cup p(K_{m}^{bot}, K_{m}^{top}) \cup p(K_{\ell}^{top}, K_{\ell}^{bot}).$$ 

Note that
\begin{equation}\label{eq:2P1-T1}
2|P_3|- |T_1|
 \geq (m-1)\cdot \min\{x_\ell,m\} 
+ (\ell-1) \cdot \min\{x_m - m,\ell\}  - \ell m - (x_\ell-\ell)(x_m-m).
\end{equation}

\begin{case}
$x_\ell \leq m$.
\end{case}

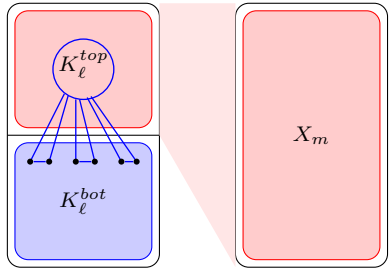
\begin{figure}[!h]
  \centering
  \begin{tikzpicture}[scale=0.5]
    \draw[rounded corners=5pt] (0,0) rectangle (4,7);
    \draw (0,3.5) -- (4,3.5);
    \draw[rounded corners=5pt] (6,0) rectangle (10,7);

    \draw[rounded corners=5pt, red, fill=red!20] (0.2,3.7) rectangle (3.8,6.8);
    \node[font=\footnotesize] at (2,5.38) {$K_{\ell}^{top}$};
    \draw[blue, line width=0.5pt] (2,5.25) circle (0.8);
    \draw[rounded corners=5pt, blue, fill=blue!20] (0.2,0.2) rectangle (3.8,3.3);
    \node[font=\footnotesize] at (2,1.75) {$K_{\ell}^{bot}$};
     
    \draw[rounded corners=5pt,red,fill=red!20] (6.2,0.2) rectangle (9.8,6.8);

    \draw[blue, line width=0.5pt] (1.5,4.6) -- (0.6,2.8);
    \draw[blue, line width=0.5pt] (1.6,4.55) -- (1.1,2.8);
    \draw[blue, line width=0.5pt] (1.8,4.45) -- (1.8,2.8);
    \draw[blue, line width=0.5pt] (1.9,4.45) -- (2.3,2.8);
    \draw[blue, line width=0.5pt] (2.1,4.5) -- (3,2.8);
    \draw[blue, line width=0.5pt] (2.2,4.53) -- (3.4,2.8);
    \node[circle, fill=black, inner sep=0.8pt] (R1) at (0.6,2.8) {};
    \node[circle, fill=black, inner sep=0.8pt] (S1) at (1.1,2.8) {};
    \node[circle, fill=black, inner sep=0.8pt] (R2) at (1.8,2.8) {};
    \node[circle, fill=black, inner sep=0.8pt] (S2) at (2.3,2.8) {};
    \node[circle, fill=black, inner sep=0.8pt] (R3) at (3,2.8) {};
    \node[circle, fill=black, inner sep=0.8pt] (S3) at (3.4,2.8) {};
  
    \draw[blue, line width=0.5pt] (R1) -- (S1);
    \draw[blue, line width=0.5pt] (R2) -- (S2);
    \draw[blue, line width=0.5pt] (R3) -- (S3);
    
     \node[font=\footnotesize] at (8,3.5) {$X_m$};
     
    \begin{scope}
      \clip (4,7) -- (4,3.5) -- (6,0) -- (6,7);
      \fill[red!10] (0,0) rectangle (10,7);
    \end{scope}
  \end{tikzpicture}
  \caption{Packing $P_4 := p(K_{\ell}^{top}\cup X_{m}) \cup p(K_{\ell}^{top},K_{\ell}^{bot}).$ Note that we use $P_4$ when $|X_m| = 2m.$ }
  \label{P4}
\end{figure}

In this case, by \eqref{eq:2P1-T1},
\begin{eqnarray}
2|P_3|- |T_1|
 &\geq&  (m-1)\cdot x_\ell 
+ (\ell-1) \cdot \min\{x_m - m,\ell\}  - \ell m - (x_\ell-\ell)(x_m-m) \nonumber \\ 
&\geq& (x_\ell - \ell)\cdot(2m-x_m) + (\ell-1)\cdot \min\{x_m-m,\ell\} -x_\ell. \label{eq:2P1-T11}
\end{eqnarray}
First, suppose that $\min\{x_m-m,\ell\} = \ell $. If $x_m <2m$, then 
\begin{eqnarray*}
2|P_3|- |T_1|
 &\geq&  (x_\ell - \ell)\cdot(2m-x_m) + (\ell-1)\cdot \ell -x_\ell \\
 &\geq& (x_\ell - \ell) + (\ell-1)\cdot \ell -x_\ell \\
 &\geq& (\ell-2)\cdot \ell \\ 
 &\geq& 0.
\end{eqnarray*}
So, if $x_m=2m$, then, by \eqref{eq:2P1-T11}, $2|P_3|- |T_1| \geq (\ell-1)\cdot \ell -x_\ell \geq (\ell-3)\cdot \ell$, which is nonnegative when $\ell \geq 3$. Hence, we may assume that $\ell=2$. Also, we may assume that $m \geq 3$ due to Theorem 1.2 of \cite{Puleo15}.  So, $2|P_3|- |T_1| \geq 2-x_{\ell}$, which is nonnegative when $x_{\ell}=2$. For $x_{\ell} \geq 3$, we define a packing $P_4$ as follows~(Figure~\ref{P4}).
$$P_4 := p(K_{\ell}^{top}\cup X_{m}) \cup p(K_{\ell}^{top},K_{\ell}^{bot}).$$

Note that, by Proposition 2, since $2m+2 \neq 5$ we have
\begin{eqnarray}
2|P_4|- |T_1|
 &\geq& \frac{2}{3}(\binom{2m+2}{2} - m-2) -m^{2}-m(x_{\ell} - 1) \nonumber \\
 &=& \frac{1}{3}(m^{2}-2-m(3x_{\ell}-7)). \label{eq:2P1-T111}
\end{eqnarray}
If $x_{\ell}=3$, then, by \eqref{eq:2P1-T111}, we have that $2|P_4| - |T_1| \geq \frac{1}{3}(m^{2}-2m-2) \geq 1/3$. For $x_{\ell}=4$, we have, by \eqref{eq:2P1-T111}, that $2|P_4| - |T_1| \geq \frac{1}{3}(m^{2}-5m-2)$, which is at least -2/3 for $m \geq 5$. Hence, we may assume that $m \leq 4$. If $m=3$, then, we have that $2|P_2| - |T_1| \geq 0$. For $m=4$, we will choose a packing $P_5$ as follows.
$$P_5 := p(K_{\ell}^{top}\cup X_{m}).$$
Then, $|P_5| = \frac{1}{3}(\binom{10}{2}- \frac{10}{2} - 1) = 13$. Note that there exist at least one edge in $K_{\ell}^{top}\cup X_{m}$ which is not selected in any triangle of $P_5$. Without loss of generality, we may assume that $c_1c_2$ is such edge. Let us define a triangle packing $P_{5}^{'}$ as follows.
$$P_{5}^{'} := P_5 \cup \{c_1c_2c_3,c_3c_4d_8\}.$$
Thus, $2|P_5^{'}| -|T_1| \geq 0$.

Now suppose that $\min\{x_m-m,\ell\} = x_m - m $. In that case,
\begin{eqnarray*}
2|P_3|- |T_1|
 &\geq&  (x_\ell - \ell)\cdot(2m-x_m) + (\ell-1)\cdot(x_m-m) -x_\ell.
\end{eqnarray*}
Note that, if $x_\ell - \ell \geq 1$, then, as $\ell \geq 2$,
\begin{eqnarray*}
2|P_3|- |T_1|
 &\geq&  2m-x_m + (x_m-m) -x_\ell \\
 &\geq& m - x_\ell \\
 &\geq& 0.
\end{eqnarray*}
Hence, we may assume that $x_\ell = \ell$.

First, suppose that $\ell + m = 5$. Since $m \geq \ell \geq 2$, then $m=3$ and $\ell=2$. So, we have that $|T_1| = 14$. For this case, we will choose a packing $P_{6}$ defined as follows. $$P_6:=p(K_{\ell})\cup p(K_{m}).$$
Note that there exists at least one edge in $K_{\ell}$ which is not selected in any triangle of $P_{6}$. Without loss of generality, we may assume that $c_1c_2$ is such edge. Analogously, we may assume that $d_5d_6$ does not lie in any triangle of $P_{6}$. 
Hence, $P_{6} \cup \{c_1c_2d_4, c_1d_5d_6\}$
has size 7. Since $|T_1|=14$, the statement follows.

\begin{figure}[!h]
  \centering
  \begin{tikzpicture}[scale=0.5]
    \draw[rounded corners=5pt] (0,0) rectangle (4,7);
    \draw (0,3.5) -- (4,3.5);
    \draw[rounded corners=5pt] (6,0) rectangle (10,7);
    \draw (6,3.5) -- (10,3.5);

    \draw[rounded corners=5pt, red, fill=red!20] (0.2,3.7) rectangle (3.8,6.8);
    \node[font=\footnotesize] at (2,5.38) {$X_{\ell}$};
    \draw[blue, line width=0.5pt] (2,5.25) circle (0.8);
    \draw[rounded corners=5pt, blue, fill=blue!20] (0.2,0.2) rectangle (3.8,3.3);
    \node[font=\footnotesize] at (2,1.75) {$K_{\ell}^{bot}$};
     
    \draw[rounded corners=5pt,red,fill=red!20] (6.2,0.2) rectangle (9.8,3.3);

    \draw[blue, line width=0.5pt] (1.5,4.6) -- (0.6,2.8);
    \draw[blue, line width=0.5pt] (1.6,4.55) -- (1.1,2.8);
    \draw[blue, line width=0.5pt] (1.8,4.45) -- (1.8,2.8);
    \draw[blue, line width=0.5pt] (1.9,4.45) -- (2.3,2.8);
    \draw[blue, line width=0.5pt] (2.1,4.5) -- (3,2.8);
    \draw[blue, line width=0.5pt] (2.2,4.53) -- (3.4,2.8);
    \node[circle, fill=black, inner sep=0.8pt] (R1) at (0.6,2.8) {};
    \node[circle, fill=black, inner sep=0.8pt] (S1) at (1.1,2.8) {};
    \node[circle, fill=black, inner sep=0.8pt] (R2) at (1.8,2.8) {};
    \node[circle, fill=black, inner sep=0.8pt] (S2) at (2.3,2.8) {};
    \node[circle, fill=black, inner sep=0.8pt] (R3) at (3,2.8) {};
    \node[circle, fill=black, inner sep=0.8pt] (S3) at (3.4,2.8) {};
  
    \draw[blue, line width=0.5pt] (R1) -- (S1);
    \draw[blue, line width=0.5pt] (R2) -- (S2);
    \draw[blue, line width=0.5pt] (R3) -- (S3);

    \begin{scope}
      \clip (4,7) -- (4,3.5) -- (6,0) -- (6,3.5);
      \fill[red!10] (0,0) rectangle (10,7);
    \end{scope}
    
    \draw[blue, line width=0.5pt] (8,1.5) circle (0.8);
    \node[font=\footnotesize][black] at (8,1.5) {$K_m^{bot}$};
    
    \draw[rounded corners=5pt, blue, fill=blue!20] (6.2,3.7) rectangle (9.8,6.8);
    \node[font=\footnotesize] at (8,5.25) {$K_m^{top}$};    
    
    \draw[blue, line width=0.5pt] (7.6,2.2) -- (6.6,4.2);
    \draw[blue, line width=0.5pt] (7.7,2.25) -- (7.1,4.2);
    \draw[blue, line width=0.5pt] (8,2.3) -- (7.8,4.2);
    \draw[blue, line width=0.5pt] (8.1,2.25) -- (8.3,4.2);
    \draw[blue, line width=0.5pt] (8.3,2.2) -- (9,4.2);
    \draw[blue, line width=0.5pt] (8.4,2.15) -- (9.5,4.2);
    
    \node[circle, fill=black, inner sep=0.8pt] (P1) at (6.6,4.2) {};
    \node[circle, fill=black, inner sep=0.8pt] (Q1) at (7.1,4.2) {};
    \node[circle, fill=black, inner sep=0.8pt] (P2) at (7.8,4.2) {};
    \node[circle, fill=black, inner sep=0.8pt] (Q2) at (8.3,4.2) {};
    \node[circle, fill=black, inner sep=0.8pt] (P3) at (9,4.2) {};
    \node[circle, fill=black, inner sep=0.8pt] (Q3) at (9.5,4.2) {};
    
    \draw[blue, line width=0.5pt] (P1) -- (Q1);
    \draw[blue, line width=0.5pt] (P2) -- (Q2);
    \draw[blue, line width=0.5pt] (P3) -- (Q3);
    
  \end{tikzpicture}
  \caption{Packing $P_7 := p(X_\ell \cup K_{m}^{bot}) \cup p( K_{\ell}^{top}, K_{\ell}^{bot}) \cup  p(K_{m}^{bot},K_{m}^{top}).$ Note that we use $P_7$ when $X_{\ell}$ = $K_{\ell}^{top}.$}
  \label{P7}
\end{figure}
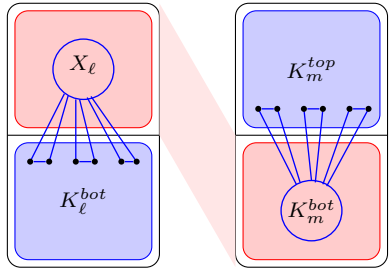

Hence, we may assume that $\ell + m \neq 5$.
We define a triangle packing $P_7$ as follows (Figure \ref{P7}).
   $$P_7 := p(X_\ell \cup K_{m}^{bot}) \cup p( K_{\ell}^{top}, K_{\ell}^{bot}) \cup  p(K_{m}^{bot},K_{m}^{top}).$$
We have
\begin{eqnarray}
2|P_7|- |T_1|
 &\geq& \frac{2}{3}(\binom{\ell + m}{2} - \frac{\ell+m}{2} - 1) -\ell m \nonumber \\
 &=& \frac{1}{3}(\ell^{2} + m^{2} - \ell m - 2(\ell + m) - 2) \nonumber \\
 &=& \frac{1}{3}((m-\ell)^{2} + (m-2)(\ell - 2) - 6) \label{eq:2P2-T1geq13m-ll+m-1l-2-6}.
\end{eqnarray}

First, suppose that $m-\ell \geq 2$, then, by \eqref{eq:2P2-T1geq13m-ll+m-1l-2-6}, we have that $2|P_7|- |T_1| \geq -\frac{2}{3}$.
Suppose now that $m-\ell =1$. Then, $m+\ell = 2\ell + 1$, which is odd. Then,  
\begin{eqnarray*}
2|P_7|- |T_1|
 &\geq& \frac{2}{3}(\binom{2\ell + 1}{2} - 4) - \ell\cdot(\ell + 1) \\
 &\geq& \frac{1}{3}(\ell^{2} - \ell - 8 ).
\end{eqnarray*}
Since $m + \ell \neq 5 $, then $\ell \geq 3$. So~$2|P_7|- |T_1| \geq -\frac{2}{3}$.

Finally, if $m-\ell =0$, then, when $\ell = 2$, by Theorem 1.2 of \cite{Puleo15}, the statement follows. So, we may assume that $\ell \geq 3$. 
Note that, by \eqref{eq:2P2-T1geq13m-ll+m-1l-2-6}, $2|P_7|- |T_1| \geq -\frac{2}{3}$ when $\ell \geq 4$. Hence, we may assume that $\ell = 3$.
In that situation, we have that $\ell = m = x_\ell = 3$.
Now, if $x_m=3$, we have that $|P_7| = \binom{3}{2} + \binom{3}{2} + \frac{1}{3}(\binom{6}{2} - 3) = 10$. Now, in $T_1$, note that we can delete the edge $c_3d_4$ because all triangles in which $c_3d_4$ lies have their third vertex in either $K_{\ell}^{top}$ or $K_m^{bot}$. Thus, $|T_1 \setminus \{c_3d_4\}| \leq 2\binom{3}{2} + 2\binom{3}{2} + 9 - 1 = 20$ and $2|P_7|- |T_1 \setminus \{c_3d_4\}| \geq 0$.
Hence, we may assume that $x_m \geq 4$. 
We define a triangle packing $P_8$ as follows.
$$P_8 := p(X_\ell \cup K_{m}^{bot} \cup \{d_3\}) \cup p( K_{\ell}^{top}, K_{\ell}^{bot}) \cup  p(K_{m}^{bot},K_{m}^{top} \setminus \{d_3\}).$$

Then, $|P_8| = \frac{1}{3}\binom{7}{2} + \binom{3}{2} + \binom{2}{2} = 11$, so $2|P_8| - |T_1| \geq 1$. 

 \begin{case} $x_{\ell}>m$.  \\
 
     We divide this case in two subcases, depending
on the relationship between $x_m-m$ and $\ell$.

 \begin{subcase} $x_m \leq m+\ell$.
 \end{subcase}
 In this case,  $\min\{x_m - m,\ell\}=x_m - m$, so, by \eqref{eq:2P1-T1},
\begin{eqnarray}
2|P_3|- |T_1|
 &\geq& (m-1)\cdot m 
+ (\ell-1) \cdot (x_m - m)  - \ell m - (x_\ell-\ell)(x_m-m)\nonumber \\
&=& 
m\cdot (m-\ell-1) 
+ (x_m-m)\cdot(2\ell - 1 -x_\ell). \label{eq:2P1-T1whenxl>=m}
\end{eqnarray}
Now, if $m-\ell \geq 2$, then, by \eqref{eq:2P1-T1whenxl>=m},
$
2|P_3|- |T_1|
 \geq m + (x_m-m)(2\ell - 1 -x_\ell)
 \geq m + m - x_m 
\geq 0.
$
Hence, we may assume that $m-\ell \leq 1$. 

If $m-\ell=1$, then, by \eqref{eq:2P1-T1whenxl>=m},
$2|P_3|- |T_1| \geq (x_m-m)\cdot(2\ell - 1-x_\ell)$,
which is nonnegative if $x_\ell < 2\ell$. So, let us assume that $x_\ell = 2\ell$.
Note that, since $3\ell + 1 \neq 5$,
\begin{eqnarray}
2|P_2|- |T_1|
 &\geq& \frac{2}{3}(\binom{3\ell + 1}{2} - \frac{3\ell + 1}{2} - 1) - \ell \cdot (\ell + 1) - 2\binom{\ell}{2} - (x_m-m)\cdot \ell \nonumber \\
 &=& \ell^{2} - 1 - \ell \cdot (x_m-m) \label{eq:2P4-T1geqll-1-xm-m}.
\end{eqnarray}
If $x_m - m \leq \ell - 1$, then, by \eqref{eq:2P4-T1geqll-1-xm-m}, we have that $2|P_2|- |T_1| \geq 0$. As $x_m - m \leq \ell $, we can assume that $x_m = m + \ell = 2m - 1$.
	
Note that if the total number of edges between
$X_{\ell} \setminus K_\ell^{top}$ and $X_{m} \setminus K_{m}^{bot}$ is less than $(x_m-m)\cdot(x_\ell - \ell)$, then $\lvert T_1 \rvert \leq \ell m + (x_\ell  - \ell)(x_m - m) - 1 +
2\binom{m}{2}+ 2\binom{\ell}{2}$. So,
\begin{eqnarray*}
2|P_2|- |T_1|
 &\geq& \ell^{2} - \ell \cdot (x_m-m) \\
 &=& 0.
\end{eqnarray*}
Hence, we may assume that $X_{\ell} \setminus K_\ell^{top}$ and $X_{m} \setminus K_{m}^{bot}$ are complete to each other. We define a triangle packing $P_{9}$ as follows.
$$P_9:= p(X_{\ell} \cup X_{m}).$$
We have
\begin{eqnarray}
2|P_9|- |T_1|
 &\geq& \frac{2}{3}(\binom{4\ell + 1}{2} - 4) -4\ell^{2}-\ell \nonumber \\
 &=& \frac{1}{3}(4\ell^{2} + \ell - 8) \nonumber \\
 &\geq& \frac{10}{3}. \nonumber
\end{eqnarray}

Now, if $m-\ell = 0$, let us assume that $\ell$ is odd because the case when $\ell$ is even follows by Theorem~2~of~\cite{Bonamy2021}. Thus $\ell \geq 3$. Now, we define a triangle packing $P_{10}$ as follows (Figure \ref{P10}).
$$P_{10}:= p(K_{m}^{bot},K_{m}^{top}) \cup p(K_{\ell}^{top}, K_{\ell}^{bot}) \cup p(K_{m}^{bot},K_{\ell}^{top}) \cup p(X_{\ell} \setminus K_{\ell}^{top}, K_{m}^{bot}).$$

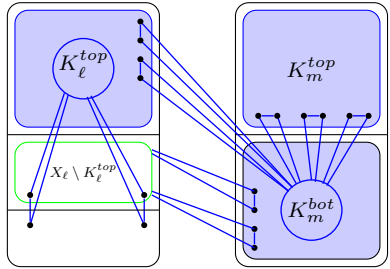
\begin{figure}[!h]
  \centering
  \begin{tikzpicture}[scale=0.5]
    \draw[rounded corners=5pt] (0,0) rectangle (4,7);
    \draw (0,3.5) -- (4,3.5);
    \draw[rounded corners=5pt] (6,0) rectangle (10,7);
    \draw (6,3.5) -- (10,3.5);
    
    \draw[rounded corners=5pt, blue, fill=blue!20] (0.2,3.7) rectangle (3.8,6.8);
    \node[font=\footnotesize] at (2,5.38) {$K_{\ell}^{top}$};
    \draw[blue, line width=0.5pt] (2,5.25) circle (0.8);
    \draw[rounded corners=5pt, green] (0.2,1.7) rectangle (3.8,3.3);
    \node[font=\footnotesize] at (2,2.5) {\resizebox{!}{0.18cm}{$X_{\ell} \setminus K_{\ell}^{top}$}};

    
    \draw[rounded corners=5pt, blue, fill=blue!20] (6.2,3.7) rectangle (9.8,6.8);
    \node[font=\footnotesize] at (8,5.25) {$K_m^{top}$}; 
    \draw[rounded corners=5pt, fill=blue!20] (6.2,0.2) rectangle (9.8,3.3);
    \node at (8,2.75) {};

    \draw[blue, line width=0.5pt] (7.6,2.2) -- (6.6,4);
    \draw[blue, line width=0.5pt] (7.7,2.25) -- (7.1,4);
    \draw[blue, line width=0.5pt] (8,2.3) -- (7.8,4);
    \draw[blue, line width=0.5pt] (8.1,2.25) -- (8.3,4);
    \draw[blue, line width=0.5pt] (8.3,2.2) -- (9,4);
    \draw[blue, line width=0.5pt] (8.4,2.15) -- (9.5,4);
    
    \draw[blue, line width=0.5pt] (1.5,4.6) -- (0.6,1.9);
    \draw[blue, line width=0.5pt] (1.6,4.55) -- (0.6,1.1);
    \draw[blue, line width=0.5pt] (2.1,4.5) -- (3.6,1.1);
    \draw[blue, line width=0.5pt] (2.2,4.53) -- (3.6,1.9);

    \draw[blue, line width=0.5pt] (7.5,2.1) -- (3.5,6.5);
    \draw[blue, line width=0.5pt] (7.5,2.1) -- (3.5,6.0);
    \draw[blue, line width=0.5pt] (7.4,2) -- (3.5,5.5);
    \draw[blue, line width=0.5pt] (7.4,2) -- (3.5,5.0);

    \node[circle, fill=black, inner sep=0.8pt] (P1) at (6.6,4) {};
    \node[circle, fill=black, inner sep=0.8pt] (Q1) at (7.1,4) {};
    \node[circle, fill=black, inner sep=0.8pt] (P2) at (7.8,4) {};
    \node[circle, fill=black, inner sep=0.8pt] (Q2) at (8.3,4) {};
    \node[circle, fill=black, inner sep=0.8pt] (P3) at (9,4) {};
    \node[circle, fill=black, inner sep=0.8pt] (Q3) at (9.5,4) {};
    
    \node[circle, fill=black, inner sep=0.8pt] (R1) at (0.6,1.9) {};
    \node[circle, fill=black, inner sep=0.8pt] (S1) at (0.6,1.1) {};
    \node[circle, fill=black, inner sep=0.8pt] (R3) at (3.6,1.1) {};
    \node[circle, fill=black, inner sep=0.8pt] (S3) at (3.6,1.9) {};
    
    \node[circle, fill=black, inner sep=0.8pt] (T1) at (3.5,6.5) {};
    \node[circle, fill=black, inner sep=0.8pt] (U1) at (3.5,6.0) {};
    \node[circle, fill=black, inner sep=0.8pt] (T2) at (3.5,5.5) {};
    \node[circle, fill=black, inner sep=0.8pt] (U2) at (3.5,5.0) {};
    \node[circle, fill=black, inner sep=0.8pt] (V2) at (6.5,2.0) {};
    \node[circle, fill=black, inner sep=0.8pt] (W2) at (6.5,1.5) {};
    \node[circle, fill=black, inner sep=0.8pt] (V3) at (6.5,1.0) {};
    \node[circle, fill=black, inner sep=0.8pt] (W3) at (6.5,0.5) {};     
    
    \draw[blue, line width=0.5pt] (P1) -- (Q1);
    \draw[blue, line width=0.5pt] (P2) -- (Q2);
    \draw[blue, line width=0.5pt] (P3) -- (Q3);
    
    \draw[blue, line width=0.5pt] (R1) -- (S1);
    \draw[blue, line width=0.5pt] (R3) -- (S3);
    
    \draw[blue, line width=0.5pt] (T1) -- (U1);
    \draw[blue, line width=0.5pt] (T2) -- (U2);
    
    \draw[blue, line width=0.5pt] (V2) -- (W2);
    \draw[blue, line width=0.5pt] (V3) -- (W3);
    
    \draw[blue, line width=0.5pt] (8,1.5) circle (0.8);
    \node[font=\footnotesize][black] at (8,1.5) {$K_m^{bot}$};
    \node at (8,4.85) {};
    
    \draw (0,1.5) -- (4,1.5);
    
    \draw[blue, line width=0.5pt] (3.8,3.1) -- (6.5,2.0);
    \draw[blue, line width=0.5pt] (3.8,3.0) -- (6.5,1.5);
    \draw[blue, line width=0.5pt] (3.8,2.0) -- (6.5,1.0);
    \draw[blue, line width=0.5pt] (3.8,1.9) -- (6.5,0.5);
  \end{tikzpicture}
  \caption{Packing $P_{10}:= p(K_{m}^{bot},K_{m}^{top}) \cup p(K_{\ell}^{top}, K_{\ell}^{bot}) \cup p(K_{m}^{bot},K_{\ell}^{top}) \cup p(X_{\ell} \setminus K_{\ell}^{top}, K_{m}^{bot}).$ }
  \label{P10}
\end{figure}

Also, since $\ell$ is odd, there exists a perfect matching between $K_{m}^{bot}$ and $K_{\ell}^{top}$ that was not used in any triangle of $P_{10}$~\cite[Lemma 4]{Bonamy2021}. The same happens between $K_{m}^{bot}$ and $K_{m}^{top}$. As the edges between $K_{\ell}^{top}$ and $X_{m} \setminus K_{m}^{bot}$ were not used, we can obtain $x_m - m = x_m - \ell$ more triangles for $P_{10}$. Let $P'_{10}$ be the resulting packing.
We have
\begin{eqnarray}
2|P'_{10}|- |T_1|
 &\geq& 2\binom{\ell}{2} + (\ell - 1)\cdot \min\{x_\ell - \ell,\ell \} + 2(x_m - \ell) - \ell^{2} - (x_\ell - \ell)\cdot (x_m - \ell) \nonumber \\
 &=& (\ell - 1)\cdot(x_\ell - \ell) + 2x_m - 3\ell - (x_\ell - \ell)(x_m - \ell) \nonumber \\
 &=& (x_\ell - \ell - 1)\cdot (2\ell - 1 - x_m) + x_m - \ell - 1.
 \label{eq:2P6-T1}
\end{eqnarray} 
If $x_\ell > \ell + 1$, then, by \eqref{eq:2P6-T1}, $2|P'_{10}|- |T_1| \geq \ell -2 \geq 0$. Since $x_\ell > m$, we may assume that $x_\ell =\ell + 1$. In that case, by \eqref{eq:2P6-T1}, if $x_m > \ell$ then 
$2|P_{10}|- |T_1| \geq x_m - \ell - 1 \geq 0$.
So, we can assume that $x_m=\ell$. We define a triangle packing $P_{11}$ as follows (Figure \ref{P11}).
$$P_{11}:= p(X_\ell \cup K_{m}^{bot}) \cup p(X_{\ell}, K_{\ell} \setminus X_{\ell}) \cup p(K_{m}^{bot},K_{m}^{top}).$$
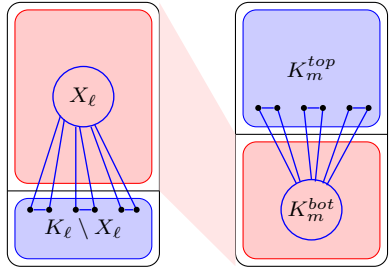
\begin{figure}[!h]
  \centering
  \begin{tikzpicture}[scale=0.5]
    \draw[rounded corners=5pt] (0,0) rectangle (4,7);
    
    \draw[rounded corners=5pt] (6,0) rectangle (10,7);
    \draw (6,3.5) -- (10,3.5);

    \draw[rounded corners=5pt, red, fill=red!20] (0.2,2.2) rectangle (3.8,6.8);
    \node[font=\footnotesize] at (2,4.5) {$X_{\ell}$};
    \draw[blue, line width=0.5pt] (2,4.5) circle (0.8);
    \draw[rounded corners=5pt, blue, fill=blue!20] (0.2,0.2) rectangle (3.8,1.8);
    \node[font=\footnotesize] at (2,1) {$K_{\ell} \setminus X_{\ell}$};
     
    \draw[rounded corners=5pt,red,fill=red!20] (6.2,0.2) rectangle (9.8,3.3);

    \draw[blue, line width=0.5pt] (1.4,4) -- (0.6,1.5);
    \draw[blue, line width=0.5pt] (1.5,3.9) -- (1.1,1.5);
    \draw[blue, line width=0.5pt] (1.8,3.7) -- (1.8,1.5);
    \draw[blue, line width=0.5pt] (1.9,3.7) -- (2.3,1.5);
    \draw[blue, line width=0.5pt] (2.2,3.7) -- (3,1.5);
    \draw[blue, line width=0.5pt] (2.3,3.75) -- (3.4,1.5);
    \node[circle, fill=black, inner sep=0.8pt] (R1) at (0.6,1.5) {};
    \node[circle, fill=black, inner sep=0.8pt] (S1) at (1.1,1.5) {};
    \node[circle, fill=black, inner sep=0.8pt] (R2) at (1.8,1.5) {};
    \node[circle, fill=black, inner sep=0.8pt] (S2) at (2.3,1.5) {};
    \node[circle, fill=black, inner sep=0.8pt] (R3) at (3,1.5) {};
    \node[circle, fill=black, inner sep=0.8pt] (S3) at (3.4,1.5) {};
  
    \draw[blue, line width=0.5pt] (R1) -- (S1);
    \draw[blue, line width=0.5pt] (R2) -- (S2);
    \draw[blue, line width=0.5pt] (R3) -- (S3);

    \begin{scope}
      \clip (4,7) -- (4,2) -- (6,0) -- (6,3.5);
      \fill[red!10] (0,0) rectangle (10,7);
    \end{scope}
    
    \draw[blue, line width=0.5pt] (8,1.5) circle (0.8);
    \node[font=\footnotesize][black] at (8,1.5) {$K_m^{bot}$};
    
    \draw[rounded corners=5pt, blue, fill=blue!20] (6.2,3.7) rectangle (9.8,6.8);
    \node[font=\footnotesize] at (8,5.25) {$K_m^{top}$};    
    
    \draw[blue, line width=0.5pt] (7.6,2.2) -- (6.6,4.2);
    \draw[blue, line width=0.5pt] (7.7,2.25) -- (7.1,4.2);
    \draw[blue, line width=0.5pt] (8,2.3) -- (7.8,4.2);
    \draw[blue, line width=0.5pt] (8.1,2.25) -- (8.3,4.2);
    \draw[blue, line width=0.5pt] (8.3,2.2) -- (9,4.2);
    \draw[blue, line width=0.5pt] (8.4,2.15) -- (9.5,4.2);
    
    \node[circle, fill=black, inner sep=0.8pt] (P1) at (6.6,4.2) {};
    \node[circle, fill=black, inner sep=0.8pt] (Q1) at (7.1,4.2) {};
    \node[circle, fill=black, inner sep=0.8pt] (P2) at (7.8,4.2) {};
    \node[circle, fill=black, inner sep=0.8pt] (Q2) at (8.3,4.2) {};
    \node[circle, fill=black, inner sep=0.8pt] (P3) at (9,4.2) {};
    \node[circle, fill=black, inner sep=0.8pt] (Q3) at (9.5,4.2) {};
    
    \draw[blue, line width=0.5pt] (P1) -- (Q1);
    \draw[blue, line width=0.5pt] (P2) -- (Q2);
    \draw[blue, line width=0.5pt] (P3) -- (Q3);

    \draw (0,2) -- (4,2);
  \end{tikzpicture}
  \caption{Packing $P_{11}:= p(X_\ell \cup K_{m}^{bot}) \cup p(X_{\ell}, K_{\ell} \setminus X_{\ell}) \cup p(K_{m}^{bot},K_{m}^{top}).$}
  \label{P11}
\end{figure}

We have
\begin{eqnarray}
2|P_{11}|- |T_1|
 &\geq& \frac{2}{3}(\binom{2\ell + 1}{2} - 4) + (\ell - 1)\cdot (\ell - 2) - \ell(\ell - 1) - \ell ^{2} \nonumber \\
 &=& \frac{1}{3}(\ell^{2} - 4\ell - 2). \label{eq:2P7-T1}
\end{eqnarray}
Note that, by \eqref{eq:2P7-T1}, if $\ell \geq 5$, then $2|P_{11}|- |T_1| \geq 1$.
Hence, we may assume that $\ell = 3$. Since $2\ell + 1=7$, we have
$2|P_{11}|- |T_1|
 \geq \frac{2}{3}\binom{2\ell + 1}{2} + (\ell - 1)(\ell - 2) - \ell(\ell - 1) - \ell ^{2} = 1$,
 and the proof follows.

\begin{subcase} $x_m > m + \ell$.
\end{subcase}

First, note that $m > \ell$ because, if $m=\ell$, then $x_m > 2m$, which is a contradiction. Now, let us define a triangle packing $P_{12}$ as follows.
$$P_{12}:= p(K_{m}^{bot},X_{\ell}) \cup p(K_{\ell}^{top},X_{m} \setminus K_{m}^{bot}) \cup p(K_{m}^{top},K_{m}^{bot}).$$ Then
\begin{eqnarray}
2|P_{12}|- |T_1|
&\geq&(x_\ell-1)\cdot \min \{m,x_\ell\}+ 
(x_m-m-1)\cdot \min\{\ell,x_m-m\} \nonumber \\
&&  - \ell m - (x_\ell-\ell)(x_m-m)-\ell(\ell-1) \nonumber \\
 &=& (x_\ell - 1)\cdot m+ 
(x_m-m-1)\cdot \ell
- \ell m - (x_\ell-\ell)(x_m-m)-\ell(\ell-1) \nonumber \\
 &=& (x_m-m)(2\ell - x_\ell) -\ell^{2} -\ell m - m + mx_\ell \nonumber \\
 &\geq& (\ell + 1)(2\ell - x_\ell) - \ell^{2} - \ell m - m + mx_\ell  \nonumber \\
 &=& \ell^{2} + 2\ell - m - \ell m + x_\ell \cdot(m-\ell-1) \nonumber \\
 &\geq& \ell^{2} + 2\ell - m - \ell m + (m+1) \cdot (m-\ell-1) \nonumber \\
  &=& (m-\ell)^2 - (m-\ell) - 1 \label{eq:2P8-T1-fin}.
\end{eqnarray}
If $m-\ell \geq 2$, then, by \eqref{eq:2P8-T1-fin}, $2|P_{12}|- |T_1| \geq 1$ and we are done. 
Hence, we may assume that $m-\ell = 1$. Since~ $x_m \geq m+\ell+1$ and $x_m \leq 2m$, we have $x_m = 2m$. 

Now, if $x_\ell =2\ell$, then,
since $x_\ell$ is even, we have that
\begin{eqnarray*}
2|P_{12}|- |T_1|
&\geq& x_\ell\cdot \min \{m,x_\ell -1\}+ 
(x_m-m-1)\cdot \min\{\ell,x_m-m\} \\
&&  - \ell m - (x_\ell-\ell)(x_m-m)-\ell(\ell-1) \\
 &=& 2\ell\cdot (\ell + 1)+ 
\ell\cdot \ell
- \ell(\ell + 1) - \ell(\ell + 1) -\ell(\ell-1) \\
 &=& \ell \\ 
 &\geq& 0.
\end{eqnarray*}
Hence, we may assume that $x_\ell \leq 2\ell - 1$. Then, we have
\begin{eqnarray*}
2|P_{4}|- |T_1|
&\geq& \frac{2}{3}(\binom{3\ell + 2}{2} - \frac{3\ell + 2}{2} - 1) - 2\binom{m}{2} - \ell m - m \cdot(x_\ell - \ell) \\
&=& \frac{3\ell^{2} - 2}{3} - (\ell + 1)(x_\ell - \ell) \\
&\geq& \frac{3\ell^{2} - 2}{3} - (\ell + 1)(\ell - 1) \\
&=& \frac{1}{3}.
\end{eqnarray*}
With that, we finish the case $x_\ell \geq \ell$.
\end{case}

\subsection{The case $x_\ell < \ell$.}

Let $K_1=K^{top}_{\ell} \cup X_{m}$ and $K_2=K^{bot}_{m} \cup X_{\ell}$. Note that $\lvert K_{1} \rvert = \ell + x_{m} $ and $\lvert K_{2} \rvert = m + x_{\ell}$. So,  without~loss of generality, let us suppose that $|K_{2}| \geq |K_{1}|$.

Let $$P_{13} := p(K_2) 
\cup  
p(X_{\ell} \cup X_m, K_{\ell}^{top}\setminus X_\ell)
\cup
p(K_{m}^{bot}, K_{m}^{top}) \cup p(K_{\ell}^{top}, K_{\ell}^{bot}).
$$

Then
$$
|P_{13}| \geq \frac{1}{3}(\binom{m+x_{\ell}}{2} - \frac{(m+x_{\ell})}{2} - \frac{3}{2}) + \frac{1}{2}(\ell-x_{\ell}-1)\cdot \min\{x_{m} + x_{\ell}, \ell-x_\ell\} + \binom{m}{2} +\binom{\ell}{2}.
$$

We choose a triangle hitting $T_2$ by taking all edges within $K_\ell^{top}$, $K_\ell^{bot}$, $K_m^{top}$, $K_m^{bot}$ as well as all edges between~$X_\ell$ and $K_m^{bot}$, and between $X_m$ and $K_\ell^{top}$. Note that 
$$
\lvert T_2 \rvert = 2\binom{m}{2} + 2\binom{\ell}{2} + mx_\ell + \ell x_m - x_\ell x_m.
$$

Hence,
\begin{equation}\label{eq:6P10-3T2}
6|P_{13}|- 3|T_2|
 \geq (m+x_{\ell})(m+x_{\ell}-2)-3 + 3(\ell-x_{\ell}-1)\cdot \min\{x_{m} + x_{\ell}, \ell-x_\ell\}  
 - 3mx_\ell - 3\ell x_m + 3x_\ell x_m.
\end{equation}
We divide the rest of the proof in two cases.


\subsubsection{The subcase $x_{m} + x_{\ell} < \ell-x_{\ell}$.}

By \eqref{eq:6P10-3T2}, we have
\begin{eqnarray}
6|P_{13}|-3|T_2|  &\geq& (m+x_{\ell})(m+x_{\ell}-2) - 3 + 3(\ell-x_{\ell}-1)\cdot (x_{m} + x_{\ell})  
 - 3mx_\ell - 3\ell x_m + 3x_\ell x_m \nonumber \\
 &=&m^2 -2m-m x_\ell - 3 -(2x_\ell^2 + 5x_\ell + 3x_m - 3\ell x_\ell). \label{eq:6P10-3T2-1}
\end{eqnarray}


Since $|K_2|\geq |K_1|$ and $x_m + x_{\ell} < \ell - x_{\ell}$, we have
\begin{equation}\label{eq:ell+x_mleqm+x_ell}
 \ell +x_m \leq m+x_\ell,
\end{equation}\\ 
and
\begin{equation}\label{eq:2x_ell+x_m+1leqm}
 2x_\ell +x_m +1 \leq \ell.
\end{equation}
By adding $2\cdot$\eqref{eq:ell+x_mleqm+x_ell}
to $1.5\cdot$\eqref{eq:2x_ell+x_m+1leqm}, we have
\begin{equation}\label{eq:2m-2-x_ell}
 2m-x_\ell \geq \frac{7}{2}x_m + \frac{\ell}{2} + \frac{3}{2}. \nonumber
\end{equation} 
Since $\ell \geq 2$, we have
\begin{equation}\label{eq:3ellgeqfrac32x_ell+fracx_m2+1}
2m-x_\ell - 2 \geq  \frac{7}{2}x_m.
\end{equation}

Using \eqref{eq:3ellgeqfrac32x_ell+fracx_m2+1} and \eqref{eq:2x_ell+x_m+1leqm} in \eqref{eq:6P10-3T2-1}, we have
\begin{eqnarray}
6|P_{13}|-3|T_2| &\geq&
m^2 -2m-m x_\ell - 3 -(2x_\ell^2 + 5x_\ell + 3x_m - 3\ell x_\ell) \label{eq:6P10-3T2-2} \\
&=&
(m-\frac{(2+x_\ell)}{2})^2 - \frac{(2+x_\ell)^2}{4}
-3- (2x_\ell^2 + 5x_\ell + 3x_m - 3\ell x_\ell) \nonumber \\
&\geq&
\frac{(7x_m)^2}{16} - \frac{(2+x_\ell)^2}{4}
-3-(2x_\ell^2 + 5x_\ell + 3x_m - 3\ell x_\ell) \nonumber \\
&=&
\frac{49x^2_m}{16}-\frac{9x^2_\ell}{4}-6x_\ell
-4-3x_m + 3\ell x_\ell \nonumber \\
&\geq&
\frac{49x^2_m}{16}-\frac{9x^2_\ell}{4}-6x_\ell
-4-3x_m + 3(2x_\ell + x_m + 1 )x_\ell \nonumber \\
&=&
(\frac{49x^2_m}{16}-3x_m)+(\frac{15x^2_\ell}{4}-3x_\ell)+
3x_mx_\ell- 4. \label{eq:6P10-3T2-3}
\end{eqnarray}
Note that, since $2|P_{13}|-|T_2|$ is integer, it suffices to prove that $6|P_{13}|-3|T_2|>-3$.
First suppose that $x_{m}\geq 2$, then, by \eqref{eq:6P10-3T2-3}, $6|P_{13}|-3|T_2|\geq 4(\frac{49}{16}) - 6 - 4 = \frac{9}{4}$. Hence, we may assume that $x_m \leq 1$. So, if $x_m=1$, then, by \eqref{eq:6P10-3T2-3} we have that $6|P_{13}|-3|T_2|> -3$ for $x_{\ell} \geq 1$. For $x_{\ell}=0$ we have, by \eqref{eq:6P10-3T2-2}, that $6|P_{13}|-3|T_2| \geq m^2-2m-6$, which is nonnegative for $m \geq 4$. Hence, we may assume that $m\leq 3$. Note that, by \eqref{eq:ell+x_mleqm+x_ell}, $m\geq \ell + 1$ and by \eqref{eq:2x_ell+x_m+1leqm}, $\ell \geq 2$. Thus, we conclude that $m \geq 3$. Since $m\leq 3$, we also conclude that $m=3$ and $\ell=2$. So, for this case we have that $2|P_{13}| - |T_2| \geq 2$. \\ \par

Finally, if $x_m=0$, then, by \eqref{eq:6P10-3T2-3}, $6|P_{13}|-3|T_2|\geq 5$ for $x_{\ell} \geq 2$. Hence, we may assume that $x_{\ell} \leq 1$. If $x_{\ell}=1$, then, by \eqref{eq:6P10-3T2-2}, $6|P_{13}|-3|T_2|\geq m^{2}-3m + 3\ell - 10$, which is nonnegative for $\ell \geq 4$ because $6|P_{13}|-3|T_2|\geq (m-2)(m-1) \geq 0$. Hence, we may assume that $\ell \leq 3$. Note that, by \eqref{eq:2x_ell+x_m+1leqm}, $\ell \geq 3$, so $\ell=3$; thus, by \eqref{eq:6P10-3T2-2}, $6|P_{13}|-3|T_2|\geq m^{2}-3m -1$, which is greater than -3 for $m\geq 3$. For $m\leq 2$, by \eqref{eq:ell+x_mleqm+x_ell}, we have that $m\geq 2$, so $m=2$. Thus, for this case, $2|P_{13}| - |T_2| \geq 2$. If $x_{\ell} =0$, then, by 
\eqref{eq:ell+x_mleqm+x_ell}, we have that $m \geq \ell$, and, by 
\eqref{eq:6P10-3T2-2}, $6|P_{13}|-3|T_2|\geq m^{2}-2m -3$, which is nonnegative for $m \geq 3$. So, let us assume that $m \leq 2$. Since $m\geq \ell$, then $\ell \leq 2$. So, by Theorem 1.2 of \cite{Puleo15}, the proof follows.


\subsubsection{The subcase $x_{m} + x_{\ell} \geq \ell-x_{\ell}$.}

By \eqref{eq:6P10-3T2}, we have
\begin{eqnarray}
6|P_{13}|-3|T_2| &\geq&
(m+x_{\ell})(m+x_{\ell}-2)-3 + 3(\ell-x_{\ell})(\ell - x_{\ell} -1) - 3(mx_\ell + \ell x_m -x_\ell x_m) \label{eq:6P10-3T2-4} \nonumber \\
&=&
m^2-2m-3+3\ell^2-3\ell+4x_\ell^2+x_\ell(1-6\ell -m)-x_m(3\ell-3x_\ell)\label{eq:6P10-3T2-5} \nonumber \\
&\geq&
m^2-2m-3+3\ell^2-3\ell+4x_\ell^2+x_\ell(1-6\ell -m)-(m+x_\ell -\ell)(3\ell-3x_\ell) \nonumber \\
&=&
m^2-2m-3+6\ell^2-3\ell-3\ell m + 7x_\ell^2-x_\ell(12\ell -2m-1) \label{eq:6P10-3T2-6} \nonumber \\
&\geq&
m^2-2m-3+6\ell^2-3\ell-3\ell m -\frac{(12\ell-2m-1)^2}{28} \nonumber \\
&=&\frac{1}{28}(24\ell^2 +24m^2-60m-60\ell-36\ell m-85).\label{eq:6|P_3|-3|T_3|-2}
\end{eqnarray}

Note that if either $m$ or $\ell$ is at least 11, then,
by \eqref{eq:6|P_3|-3|T_3|-2}, $6|P_{13}|-3|T_2|$ is greater than -3 and the proof follows.
Hence, we may assume that $m,\ell \leq 10$.

Let us define the following triangle packings.
$$P_{14} := 
p(K_{m}^{bot}, K_{m}^{top}) \cup p(K_{\ell}^{top}, K_{\ell}^{bot}) \cup p(K_{\ell}^{top} \setminus X_{\ell} \cup K_{m}^{bot} \setminus X_{m}, X_{\ell} \cup X_{m}),$$
$$P_{15}^{\ell} := p(K_{\ell}^{top}) 
\cup  
p(X_{\ell}, K_{m}^{bot})
\cup
p(K_{m}^{bot}, K_{m}^{top}) \cup p(K_{\ell}^{top}, K_{\ell}^{bot}),$$
$$P_{15}^{m} := p(K_{m}^{bot}) 
\cup  
p(X_{m}, K_{\ell}^{top})
\cup
p(K_{m}^{bot}, K_{m}^{top}) \cup p(K_{\ell}^{top}, K_{\ell}^{bot}),$$
$$P_{16}^{\ell} := p(K_{\ell}) 
\cup  
p(X_{\ell}, K_{m}^{bot})
\cup
p(K_{m}^{bot}, K_{m}^{top}),$$
$$P_{16}^{m} := p(K_{m}) 
\cup  
p(X_{m}, K_{\ell}^{top})
\cup
p(K_{\ell}^{top}, K_{\ell}^{bot}),$$
$$P_{17}^{\ell} := p(K_{\ell}) 
\cup  
p(X_{\ell}, K_{m}^{bot} \setminus d_{|K_m|})
\cup
p(K_{m}^{bot}\setminus d_{|K_m|}, K_{m}^{top} \cup d_{|K_m|}),$$ and
$$P_{17}^{m} := p(K_{m}) 
\cup  
p(X_{m}, K_{\ell}^{top} \setminus c_1)
\cup
p(K_{\ell}^{top} \setminus c_1, K_{\ell}^{bot}\cup c_1).$$

Let

$$f_1(\ell,m,x_\ell,x_m)=6|P_{13}|-3|T_2|,$$
$$f_2(\ell,m,x_\ell,x_m)=6|P_{14}|-3|T_2|,$$
$$f_3(\ell,m,x_\ell,x_m)=6|P^{\ell}_{15}|-3|T_2|,$$
$$f_4(\ell,m,x_\ell,x_m)=6|P^{m}_{15}|-3|T_2|,$$
$$f_5(\ell,m,x_\ell,x_m)=6|P_{16}^{\ell}|-3|T_2|,
$$
$$f_6(\ell,m,x_\ell,x_m)=6|P_{16}^{m}|-3|T_2|,$$
$$f_7(\ell,m,x_\ell,x_m)=6|P_{17}^{\ell}|-3|T_2|,
$$
and
$$f_8(\ell,m,x_\ell,x_m)=6|P_{17}^{m}|-3|T_2|.
$$
As $m,\ell \leq 10$; we have that $x_\ell < \ell$;
$\ell + x_m \leq m + x_{\ell}$;
$\ell-x_\ell \leq x_m + x_{\ell}$; and $x_m <m$. We can verify, by exhaustive search, that at least one of $\{f_1, \cdots ,f_8\}$ is greater than $-3$, for every
such $(\ell,m,x_\ell,x_m)$, except when $(\ell,m,x_\ell,x_m) \in \{(1,2,0,1),(2,1,1,0),(2,2,1,1),(2,3,1,2),(2,5,1,4),(3,2,2,1),(3,3,2,1),(3,4,2,3), \\ (3,6,2,5),(4,3,3,2),(5,2,4,1),(6,3,5,2)\}$.

If $(\ell,m,x_\ell,x_m) \in \{(2,5,1,4),(3,4,2,3),(3,6,2,5),(4,3,3,2),(5,2,4,1),(6,3,5,2)\}$,
then $\ell-x_{\ell}=1$ and $m-x_m=1$.
So, $K_{\ell}^{top} \cup K_{m}^{bot}$ is a clique missing a single edge, and there exists a packing with $|p(K_{\ell}^{top} \cup K_{m}^{bot})| - 1$ triangles in this clique, say $P'_{18}$.

Let $$P_{18} := P'_{18}
\cup
p(K_{m}^{bot}, K_{m}^{top}) \cup p(K_{\ell}^{top}, K_{\ell}^{bot}).
$$
We have
\begin{eqnarray*}
2|P_{18}|-|T_2| &\geq&
2|P'_{18}|- (mx_\ell + \ell x_m -x_\ell x_m),
 \label{eq:6P10-3T2-4} 
\end{eqnarray*}
which is at least 1 for $(\ell,m,x_\ell,x_m) \in \{(2,5,1,4),(3,4,2,3),(3,6,2,5),(4,3,3,2),(5,2,4,1),(6,3,5,2)\}$.
Also, note that if $(\ell,m,x_\ell,x_m) \in \{(1,2,0,1),(2,1,1,0),(2,2,1,1)\}$, then these cases were solved~by~Theorem~1.2~of \cite{Puleo15}. 
Hence, we may assume that $(\ell,m,x_\ell,x_m) \in \{(2,3,1,2),(3,2,2,1),(3,3,2,1)\}$. Since the case $(\ell,m,x_\ell,x_m) = (3,2,2,1)$ is symmetric to the case $(\ell,m,x_\ell,x_m) = (2,3,1,2)$, then, without loss of generality, we can suppose that $(\ell,m,x_\ell,x_m) \in \{(2,3,1,2),(3,3,2,1)\}$. In this direction,
we let $P^{'}_{19}:=p(K_{\ell})\cup p(K_{m})$.
Note that there exists at least one edge in $K_{\ell}$ which is not selected in any triangle of $P^{'}_{19}$. Without loss of generality, we may assume that $c_1c_2$ is such edge. Analogously, we may assume that $d_{4}d_{5}$ does not lie in any triangle of $P^{'}_{19}.$
Let us define a triangle packing $P_{19}$ as follows $$P_{19}:=P^{'}_{19} \cup \{c_1c_2d_6, c_1d_4d_5\}.$$
Then, $2|P_{19}| - |T_2| \geq 1$ for $(\ell,m,x_\ell,x_m) \in \{(2,3,1,2),(3,3,2,1)\}$, which concludes the proof of the theorem.

\section{Concluding remarks} \label{sec:remarks}

In this work, we have solved Tuza's conjecture for any co-chain graph with both sizes multiple of two. This result was obtained by analyzing the parity of variables $\ell$, $m$, $x_{\ell}$, $x_m$ defined at the beginning of Section~3 and finding adequate packing and hitting configurations that satisfy Tuza's inequality in each case. We believe this result can also be extended to any co-chain graph and moreover, for mixed unit interval graphs, a known superclass of co-chain graphs.


\begin{question}
Does Tuza's conjecture holds for any co-chain graph?

\end{question}

\begin{question}
Does Tuza's conjecture holds for any mixed unit interval graph?

\end{question}

Question 1 seems tedious to tackle with the method presented here, as the cases when both sides of the partition are not even give birth to several new cases. It seems that in order to solve Question 2, we must first solve Question 1 first. So, we believe that a solution for Tuza's conjecture is far from trivial in chordal graphs. However, we think the present techniques can be generalized to solve particular cases for other  subclasses of chordal graphs, by introducing a parity parameter in the conjecture.


\bibliographystyle{plain}
\bibliography{bibliografia}

\end{document}